\documentclass[12pt,a4paper]{article}
\usepackage{amsmath,amsxtra,amsfonts,amssymb,latexsym,amsthm,amscd,verbatim,amstext}
\usepackage[mathscr]{eucal}
\newfont{\tricyi}{wncyi10 at 12pt}
\newfont{\tricyr}{wncyr10 at 12pt}

\theoremstyle{plain}
\newtheorem{Th}{Theorem}
\newtheorem{Md}{Proposition}
\newtheorem{Lem}{Lemma}
\newtheorem{Cor}{Corollary}
\theoremstyle{definition}
\newtheorem{Def}{Definition}
\newtheorem{Note}{Remark}

\newcommand{\bpr}{\begin{proof}}
\newcommand{\epr}{\end{proof}}

\begin{document}

\centerline {{\bf  Semilinear elliptic degenerate equations with critical exponent}}

\vskip 0.7cm
\begin{center}
N. M. Tri$^1$, D. A. Tuan$^2$ 
\end{center}
\begin{center}
{\footnotesize
$^1$ Institute of Mathematics, 
Vietnam Academy of Science and Technology,\\
18 Hoang Quoc Viet, Cau Giay, Hanoi, Vietnam.\\
$^2$ University of Sciences, Vietnam National University, Hanoi\\
334 Nguyen Trai, Thanh Xuan, Hanoi, Vietnam.}
\end{center}
\vskip 0.5cm

{\bf Abstract.} In this paper we are mainly concerned with  nontrivial positive solutions to the Dirichlet problem for the degenerate elliptic equation
 \begin{gather}
 -\frac{\partial^2 u}{\partial x^2} -\left|x\right|^{2k}\frac{\partial^2 u}{\partial y^2}=|x|^{2k}u^p+f(x,y,u) \quad\text{ in }\Omega, \
 u=0  \quad\text{ on }\partial\Omega,\label{equ0}
 \end{gather}
where $\Omega$ is a bounded domain with smooth boundary in $\mathbb{R}^2, \Omega \cap \{x=0\}\ne \emptyset,$ $k\in\mathbb N,$ $f(x,y,0)=0,$ and $p=(4+5k)/k$ is the critical exponent. Recently, the equation \eqref{equ0} was investigated in \cite{Luyen:2023} for the subcritical case based on a new result obtained in \cite{NTT:2022} on embedding theorem of weighted Sobolev spaces. In the critical case considered in this paper we will essentially use the optimal functions and constants found in \cite{NTT:2022}. 

\footnotetext[3]{{\it e-mail address}: datuan1105@gmail.com}

\section{Introduction}
Semilinear boundary value problems for degenerate elliptic equations were  extensively studied for more than 30 years, see for example \cite{D. Jerison:1981a}, \cite{D. Jerison:1987a}, \cite{N. M. Tri:1998a}, \cite{N. M. Tri:1998b}, \cite{thuy:2012}, \cite{Luyen:2015}, \cite{Hua:2020b}, \cite{N. M. Tri:2010a}, \cite{N. M. Tri:2014a}, \cite{Nga:2023a}  and the references therein.  The Yamabe problem on the Heizenberg group  and CR manifolds were first treated in \cite{D. Jerison:1981a}, \cite{D. Jerison:1987a}. The 
 existence and nonexistence of nontrivial solutions to the Dirichlet boundary value problems for semilinear degenerate elliptic equations  involving the Grushin operator was first considered in \cite{N. M. Tri:1998a}, \cite{N. M. Tri:1998b}. Then the results of \cite{N. M. Tri:1998a}, \cite{N. M. Tri:1998b} were then generalized for semilinear  equations containing the strongly degenerate elliptic operator in \cite{thuy:2012}. The non-linear terms of equations in \cite{N. M. Tri:1998a}, \cite{N. M. Tri:1998b}, \cite{thuy:2012} in most interesting cases do not depend on the space variable.  Semilinear equations involving the $\Delta_\gamma$-Laplace operator was considered in \cite{Lanconelli:2012}. The class of  $\Delta_\gamma$-Laplace operators containing both the Grushin and strongly degenerate elliptic operators. The multiplicity of solutions for semilinear $\Delta_\gamma$ equations was considered in \cite{Luyen:2015}. Estimates of Dirichlet eigenvalues for degenerate $\Delta_\gamma$-Laplace operator was given in \cite{Hua:2020b}. Semilinear equations involving the Grushin operator with critical exponent was treated in \cite{Huong:2023a}. Nontrivial solutions for degenerate elliptic equations in a solid torus was investigated in a recent paper \cite{Nga:2023a}. In \cite{Luyen:2023} the authors considered semilinear degenerate elliptic equations with nonlinear terms depending on the space variables. It turns out that the critical exponents for this problem may be higher than the ones for the cases treated in \cite{N. M. Tri:1998a}, \cite{N. M. Tri:1998b}. Actually, in \cite{Luyen:2023} the authors dealt with subcritical cases based on a new results obtained in the paper \cite{NTT:2022}. Inspired by the work \cite{Brezis:1984}, in this work we consider the  critical case for the problem:
\begin{align}
	-\Delta_G u& =|x|^{2k}u^p+f(x,y,u) \quad\text{ in }\Omega, \label{equ1}\\
	u & > 0 \quad \text{ in } \Omega\setminus\{x=0\}, \label{equ2}\\
	u& =0  \quad\text{ on }\partial\Omega,\label{equ3}
\end{align} 
where $\Delta_G=\partial^2_x+|x|^{2k}\partial^2_y,$ $\Omega$ is a bounded domain with smooth boundary in $\mathbb{R}^2, \Omega \cap \{x=0\}\ne \emptyset,$ $k \in\mathbb N,$ $f(x,y,0)=0,$ and $p=(4+5k)/k$ - the critical exponent. In this critical case we will essentially use the optimal functions and constants found in \cite{NTT:2022}.

The structure of our note is as follows: In Section 2, we present our main results. They  consist of two theorems. In Section 3, we prove an embedding proposition for weighted Sobolev spaces associated with the problem. In \linebreak Section 4, we give the proof of Theorem \ref{Case1}. In Section 5, we give the proof of Theorem \ref{Case2}. Finally in Section 6, we give some results about nonexistence solution when the geometry of $\Omega$ is special. In this section, the fundamental solution of the Grushin operator $\Delta_G$ (in \cite{N. M. Tri:1998a}, \cite{N. M. Tri:2002}, \cite{Beals:1998}) is usefull.

\section{Main Results}
Let $k\in\mathbb N, p=(4+5k)/k.$ We consider the problem \eqref{equ1} - \eqref{equ3} in two cases:
\begin{itemize}
\item Case 1: $f(x, y, \xi)=\lambda|x|^{2\beta}\xi, \beta>-1/2, \lambda\in\mathbb R$;
\item Case 2: $f(x, y, \xi)=\mu|x|^{2\beta}\xi^q+|x|^{2k}h(x, y, \xi),$ $k\ge\max\{0, \beta\},$ $\beta>-1/2$, $1<q<p, \mu\in\mathbb R$;
\end{itemize}
where $f, h: \Omega\times[0, \infty)\to \mathbb R, h(x, y, 0)=0.$

For stating our main results we need some function spaces as follows.\\
\begin{Def} For $S^2_{1, 0}(\Omega)$ is the completion of $C^1_c(\Omega)$ in the norm
$$||u||_{S^2_{1,0}(\Omega)}=\left(\int_\Omega |\nabla_G u|^2dxdy\right)^{1/2}$$
where $\nabla_Gu=(u_x, |x|^k u_y).$ 
\end{Def}
\begin{Def} For $\beta\in\mathbb R, q\ge 1,$ $L^q_\beta(\Omega)$ is the space of measurable functions $u$ on $\Omega$ such that $|x|^{2\beta}|u|^q\in L^1(\Omega).$ $L^q_\beta(\Omega)$ is Banach space with the norm
$$||u||_{L^q_\beta(\Omega)}=\left(\int_\Omega|x|^{2\beta}|u|^qdxdy\right)^{1/q}.$$
For $\beta=0$ we have $L^q_0(\Omega)=L^q(\Omega).$
\end{Def}

For Case 1, we set
\begin{equation}\label{eigen}
\lambda_1=\inf\limits_{u\in S^2_{1, 0}(\Omega)\atop ||u||_{L^2_\beta(\Omega)}=1} ||u||^2_{S^2_{1, 0}(\Omega)}.
\end{equation}
\begin{Def}
$\Omega$ is called $G-$starshaped if 
$$\mathcal T\cdot\nu \ge 0 \text{ on } \partial \Omega,$$
where $\mathcal T=(x, (1+k)y)$ and $\nu=(\nu_1, \nu_2)$ is the unit outward normal vector to $\partial\Omega.$
\end{Def}
\begin{Def}
We define $S_2(\overline\Omega)$ as the linear space of functions $u\in C^1(\overline{\Omega})$
such that
\begin{align*}
\frac{\partial^2 u}{\partial x^2}, \left|x\right|^{2k}\frac{\partial^2 u}{\partial y^2}\; (\text{in distribution sense})
\end{align*}
are continuous in $\Omega$ and can be continuously extended to $\overline\Omega$.
\end{Def}
Our main result in Case 1 is as follows.
\begin{Th}\label{Case1}
Let $\lambda\in\mathbb R, \beta>-1/2$ and 
$$f(x, y, \xi)=\lambda|x|^{2\beta}\xi.$$ 
\begin{itemize}
\item[(a)] For $\lambda\ge \lambda_1$, the problem \eqref{equ1} - \eqref{equ3} has no solution in $S^2_{1, 0}(\Omega)$.
\item[(b)] For $\lambda\le  0$ and $\Omega$ is $G-$starshaped the problem \eqref{equ1} - \eqref{equ3} has no solution in $S_2(\overline{\Omega}).$
\item[(c)] For $0<\lambda<\lambda_1$ and $k\ge 2(\beta+1)$, the problem \eqref{equ1} - \eqref{equ3} has a  solution in $S^2_{1, 0}(\Omega)$.
\end{itemize}
\end{Th}
For Case 2, we assume that $h: \Omega\times[0, \infty) \to\mathbb R$ is a Caratheodory function satisfying:
\begin{itemize}
\item[(h1)] $h(x, y, 0)=0$ for  almost $(x, y)$ in $\Omega$,
\item[(h2)] the following limits hold uniformly for almost $(x, y)$ in $\Omega$
$$\lim_{\xi\to 0_+}\frac{h(x, y, \xi)}{\xi}=0, \lim_{\xi\to+\infty}\frac{h(x, y, \xi)}{\xi^p}=0.$$
\end{itemize}
We define 
\begin{equation}\label{eigen2}
S=\inf\limits_{u\in S^2_{1, 0}(\mathbb R^2)\atop ||u||_{L^{p+1}_k(\mathbb R^2)}=1} ||u||^2_{S^2_{1, 0}(\mathbb R^2)}
\end{equation}
and the following energy integral:
\begin{align}
\Psi(u) = \int_\Omega \left(\frac{|\nabla_G u|^2}{2}-\frac{|x|^{2k}|u|^{p+1}}{p+1} - \mu \frac{|x|^{2\beta}|u|^{q+1}}{q+1}-|x|^{2k}H(x, y, u)\right)dxdy
\end{align}
where $$H(x, y, \xi)=\begin{cases} \int_0^\xi h(x, y, s)ds & \text{ for } \xi\ge 0,\\ 0 & \text{ for } \xi<0.\end{cases}$$
With these definitions we can state our main result in Case 2 as follows.
\begin{Th}\label{Case2} Let $\beta>-1/2$, $1<q<p, \mu\ge 0$ and 
$$f(x, y, \xi)=\mu|x|^{2\beta}\xi^q+|x|^{2k}h(x, y, \xi).$$
Suppose that $h$ satisfies (h1)-(h2) and the embedding $S^2_{1, 0}(\Omega)\hookrightarrow L^{q+1}_\beta(\Omega)$ is compact. Then if there is a $v_0\in S^2_{1, 0}(\Omega)$ such that $v_0\ge 0, v_0\not\equiv 0,$ and
\begin{equation}\label{asmpC2}
\sup_{t\ge 0}\Psi(tv_0)<\frac{1+k}{2+3k}S^{\frac{2+3k}{2+2k}}
\end{equation}
the problem \eqref{equ1} - \eqref{equ3} has a  solution in $S^2_{1, 0}(\Omega)$.
\end{Th}
As consequence of Theorem \ref{Case2} we have the following corollary.
\begin{Cor}\label{Cor}
Let $h=0$ and $\beta>-1/2$, $1<q<p, \mu> 0$ such that
$$2k(q-1)+q<2\beta(p-1)+p.$$
 The problem \eqref{equ1} - \eqref{equ3} has a  solution in $S^2_{1, 0}(\Omega)$ if one of the following conditions  holds.
\begin{itemize}
\item[(a)] $(q+1)k>4(\beta+1).$
\item[(b)] $\mu>\mu_0$ for large enough $\mu_0>0.$
\end{itemize}
\end{Cor}
\section{Embedding Proposition}
\begin{Md}\label{imb}
\begin{itemize}
\item[(i)] For $\beta_1>\beta_2, q_1>q_2>1,$ $q_2(2\beta_1+1)< q_1(2\beta_2+1)$ or $\beta_1\le \beta_2, q_1\ge q_2\ge 1$ the following embedding
$$L^{q_1}_{\beta_1}(\Omega)\hookrightarrow L^{q_2}_{\beta_2}(\Omega)$$
is continuous.
\item[(ii)] For $\theta\in (0, 1),$ and $q_j\ge 1, \beta_j\in\mathbb R$, $j=1, 2$. Let
$$\frac{1}{q}=\frac{\theta}{q_1}+\frac{1-\theta}{q_2}, \frac{\beta}{q}=\frac{\theta\beta_1}{q_1}+\frac{(1-\theta)\beta_2}{q_2}. $$
Then the following interpolation inequality holds:
$$||u||_{L^q_\beta(\Omega)}\le ||u||^\theta_{L^{q_1}_{\beta_1}(\Omega)}||u||^{1-\theta}_{L^{q_2}_{\beta_2}(\Omega)}.$$
\item[(iii)] For $\beta>-1/2$ the following embedding
$$S^2_{1, 0}(\Omega) \hookrightarrow L^2_\beta(\Omega)$$
is compact. 
\item[(iv)] (Sobolev embedding) The following embedding
$$S^2_{1, 0}(\Omega)\hookrightarrow L^{p+1}_{k}(\Omega)$$
is continuous. 
\item[(v)] For $1<q<p$ and 
$2k(q-1)+q<2\beta(p-1)+p$
the following embedding
$$S^2_{1, 0}(\Omega)\hookrightarrow L^{q+1}_{\beta}(\Omega)$$
is compact.
\end{itemize}
\end{Md}
\begin{proof}
Note that $\Omega$ is bounded, 
$$\int_\Omega |x|^a dxdy<\infty \text{ for } a>-1,$$
 and using the Holder inequality, it is not difficult to prove (i), (ii). Using the Sobolev inequality in \cite{NTT:2022} we get the continuous embedding $S^2_{1, 0}(\Omega)\hookrightarrow L^{p+1}_k(\Omega)$ in (iv). \\
Next we will prove (iii). As in \cite{N. M. Tri:1998a} we have the embedding $S^2_{1, 0}(\Omega)
\hookrightarrow L^2(\Omega)$ is compact, so using (i) the embedding in (iii) is proved for $\beta\ge 0.$ For $\beta\in(-1/2, 0)$ there is a $\gamma\in(0, 1)$ such that $\gamma+2\beta>0.$ Note that $\Omega$ is bounded and using the Hardy inequality (\cite{Opic}) we have
\begin{align}
||u||^2_{L^2_{\beta}(\Omega)} & = \int_{|x|<\epsilon\atop (x, y)\in\Omega}|x|^{\gamma+2\beta}\frac{|u|^2}{|x|^\gamma}dxdy+ \int_{|x|>\epsilon\atop (x, y)\in\Omega}|x|^{2\beta}|u|^2dxdy\notag\\
& \le \epsilon^{\gamma+2\beta}C_{\Omega, \gamma}\int_{\Omega}|u_x|^2dxdy+ \max\{\epsilon^{2\beta}, 1\} \int_\Omega|u|^2dxdy\notag\\
& \le C_{\Omega, \gamma}\epsilon^{\gamma+2\beta}||u||^2_{S^2_{1, 0}(\Omega)}+C_\epsilon||u||^2_{L^2(\Omega)}, \forall \epsilon>0. \label{comeps}
\end{align}
Using the compact embedding $S^2_{1, 0}(\Omega)\hookrightarrow L^2(\Omega)$  and \eqref{comeps} the embedding in (iii) is proved for $-1/2<\beta<0.$\\
Next we prove (v). For $k, \beta$, $p, q$ as in (v) there is a $q_0\in (q, p)$ such that
$$2(k-\beta)(q-1)<(2\beta+1)(q_0-q).$$
So 
$$\beta_0=-\frac{k(q-1)-\beta(q_0-1)}{q_0-q}>-\frac{1}{2}.$$
Using the interpolation inequality (ii) with
$$p_1=2, p_2=q_0+1, \beta_1=\beta_0, \beta_2=k \text{ and } \theta=\frac{2(q_0-q)}{(q_0-1)(q+1)};$$
and (iii), (iv) we get (v).
\end{proof}
\begin{Note}\label{Rem1} For Sobolev embedding in (iv) we use the Sobolev inequality in \cite{NTT:2022}:
\begin{equation}\label{Sob}
S||u||_{L^{p+1}_k(\mathbb R^2)}^{2}\le ||u||_{S^2_{1, 0}(\mathbb R^2)}^2, \forall u\in S^2_{1, 0}(\mathbb R^2)
\end{equation}
where $S$ is the constant in \eqref{eigen2}, and we consider $S^2_{1, 0}(\Omega)\subset S^2_{1, 0}(\mathbb R^2)$ by setting $u=0$ outside $\Omega.$ $S$ is also the best constant of Sobolev inequality \eqref{Sob} and the equality sign holds when
\begin{equation}\label{extrem}
u(x, y)=U_\epsilon(x, y)=(\epsilon^2+r^2)^{-\frac{k}{2(k+1)}}, \epsilon>0, r=(x^{2(k+1)}+(k+1)^2y^2)^{1/2}.
\end{equation}
\end{Note}
\begin{Note}\label{Rem2} From the compact embedding (ii) and variational principle we obtain that the constant $\lambda_1$ in \eqref{eigen} is the first  eigenvalue of the eigenvalue problem:
\begin{align}
-\Delta_Gu & =\lambda |x|^{2\beta}u \text{ in } \Omega, \label{eig1}\\
u & = 0 \text{ on } \partial\Omega. \label{eig2}
\end{align}
Moreover note that $|\nabla_G|u||\le |\nabla_G u|$ and using Strong Maximal Principle, there is a $\varphi_1\in S^2_{1, 0}(\Omega), \varphi_1\ge 0$ such that $\varphi_1$ is a eigenfunction of \eqref{eig1}-\eqref{eig2} with respect to $\lambda_1$, and $\varphi_1>0$ a.e. in $\Omega\setminus\{x=0\}.$
\end{Note}
\begin{Note} Here we give an other proof for the compact embedding $S^{2}_{1,0}(\Omega)\hookrightarrow L^{r}_k(\Omega),$ $1\le r<p+1,$ in \cite{Luyen:2023}. This compact embedding and the compact embedding in Proposition \ref{imb} (iii) are the best. However we do not know whether the compact embedding in Proposition \ref{imb} (v) is the best.
\end{Note}
\section{Proof of Theorem \ref{Case1}}
Using Remark \ref{Rem2} we have
\begin{equation}\label{eig3}
\int_\Omega \nabla_G\phi \cdot \nabla_G\varphi_1dxdy=\lambda_1\int_\Omega\phi\varphi_1dxdy, \forall \phi\in S^2_{1, 0}(\Omega).
\end{equation}
So we can prove Theorem \ref{Case1} (a) as follows.
\begin{proof}[Proof of Theorem \ref{Case1} (a)] Assume that the problem \eqref{equ1}-\eqref{equ3} has a  solution $u\in S^2_{1, 0}(\Omega)$. Then we have
$$\int_\Omega \nabla_G u\cdot  \nabla_G\varphi_1dxdy=\int_\Omega |x|^{2k}u^p\varphi_1dxdy+\lambda\int_\Omega |x|^{2\beta}u\varphi_1dxdy.$$
Then using \eqref{eig3} we get
$$\int_\Omega |x|^{2k}u^p\varphi_1dxdy=(\lambda_1-\lambda)\int_\Omega |x|^{2\beta}u\varphi_1dxdy.$$
For $\lambda\ge \lambda_1$ we have a contradition since $u>0, \varphi_1>0$ a.e. in $\Omega.$\\
\end{proof}
In order to prove Theorem \ref{Case1} for $\lambda\le 0$ and $\Omega$ is $G-$starshaped we need the following Pohozaev-type indentity.
\begin{Md}\label{MdPo}
Assume that $u\in S_2(\overline{\Omega})$ is a solution of the problem \eqref{equ1}-\eqref{equ3}. Then we have
\begin{equation}\label{Pohoz}
\frac{1}{2}\int_{\partial\Omega} (\mathcal T\cdot\nu)(\nu_1^2+|x|^{2k}\nu_2^2)|\partial_\nu u|^2dS=\lambda(1+\beta)\int_\Omega |x|^{2\beta}u^2dxdy.
\end{equation}
\end{Md}
\begin{proof}
The proof of this lemma is similar to that of  Lemma 1 in \cite{thuy:2012}. We omit the details.
\end{proof}
\begin{proof}[Proof of Theorem \ref{Case1} (b)]
Note that $\Omega$ is $G-$starshaped, i.e. $\mathcal T\cdot\nu\ge 0$ on $\partial\Omega,$ using Pohozaev-type identity \eqref{Pohoz}, then as in \cite{Lanconelli:2012} the proof is done.  
\end{proof}
In order to prove Theorem \ref{Case1} (c) for $0<\lambda<\lambda_1$ we introduce some other constants:
\begin{align}
S_0:=S_0(\Omega) & = \inf_{u\in S^2_{1, 0}(\Omega)\atop ||u||_{L^{p+1}_k(\Omega)}=1}||u||^2_{S^2_{1, 0}(\Omega)}, \label{S0}\\
S_\lambda:= S_\lambda(\Omega) & = \inf_{u\in S^2_{1, 0}(\Omega)\atop ||u||_{L^{p+1}_k(\Omega)}=1}\left(||u||^2_{S^2_{1, 0}(\Omega)}-\lambda||u||^2_{L^2_\beta(\Omega)}\right).\label{Slam}
\end{align}
\begin{Note}\label{Rem3} It is obvious that $S_\lambda\le S_0.$ From \eqref{eigen}, for $\lambda<\lambda_1$ we have $S_\lambda\ge 0.$ If $V\in S^2_{1, 0}(\Omega)$ is a minimizer of \eqref{Slam}, then $|V|\in S^2_{1, 0}(\Omega)$ is also a minimizer of \eqref{Slam} since $|\nabla_G|V||\le |\nabla_G V|$. Moreover there is a constant $\tilde{\lambda}>0$ such that
$$\int_\Omega \nabla_G |V|\cdot \nabla_G\varphi dxdy - \lambda\int_\Omega |x|^{2\beta}|V|\varphi-\tilde{\lambda}\int_\Omega |x|^{2k}|V|^p\varphi dxdy=0, \forall \varphi\in S^2_{1, 0}(\Omega).$$
So $v=\tilde{\lambda}^{\frac{1}{p-1}}|V|$ is a solution of \eqref{equ1}-\eqref{equ3}. If $V\not\equiv 0$, using Maximum Principle,  $v>0$ in $\Omega\setminus\{x=0\}.$
\end{Note}
\begin{Lem} For $\lambda>0$ we have
\begin{equation}\label{import}
S_\lambda<S=S_0.
\end{equation}
\end{Lem}
\begin{proof}
Firstly, we prove that $S=S_0.$ Indeed, from $S^2_{1, 0}(\Omega)\subset S^2_{1, 0}(\mathbb R^2)$ that
\begin{equation}\label{S_0}
S \le S_0.
\end{equation}
\noindent Since $0\in \Omega$ and $\Omega$ is open, there is an $R>0$ such that 
$$B_R=\{(x, y):\; |x|^{2(k+1)}+(k+1)^2y^2<R^2\}\subset\Omega.$$
Let $\varphi\in C^\infty_c(\Omega; [0, 1])$ such that $\varphi=1$ in $B_R.$ For $\epsilon>0$ small enough, consider
$$u_\epsilon(x, y)=\varphi(x, y)\Phi_\epsilon(r), r=(|x|^{2(k+1)}+(k+1)^2y^2)^{1/2},$$
where $\Phi_\epsilon(r)=U_\epsilon(x, y)=(\epsilon^2+r^2)^{-\frac{k}{2(k+1)}}$ are extremal functions of Sobolev inequality \eqref{Sob}. By calculating we get
\begin{align*}
||u_\epsilon||^{p+1}_{L^{p+1}_k(\Omega)} & = \int_\Omega(|\varphi|^{p+1}-1)|x|^{2k}|\Phi_\epsilon(r)|^{p+1}dxdy+\int_\Omega |x|^{2k}|\Phi_\epsilon(r)|^{p+1}dxdy\\
& = O(1)+\int_{\mathbb R^2} |x|^{2k}|\Phi_\epsilon(r)|^{p+1}dxdy\\
& = O(1) +\epsilon^{-\frac{(p+1)k}{2(k+1)}}||U_1||^{p+1}_{L^{p+1}_k(\mathbb R^2)}\quad (\text{note that } p=(4+5k)/k). 
\end{align*}
So 
\begin{equation}\label{ueps_p}
||u_\epsilon||_{L^{p+1}_k(\Omega)}^2=O(1)+\epsilon^{-\frac{k}{k+1}}||U_1||_{L^{p+1}_k(\mathbb R^2)}^2.
\end{equation}
Since $\nabla_G u_\epsilon=U_\epsilon(\nabla_G\varphi)+\varphi(\nabla_G U_\epsilon)$, we have
\begin{align}
||u_\epsilon||^2_{S^2_{1, 0}(\Omega)} & = O(1) + \int_{\mathbb R^2}|\nabla_GU_\epsilon|^2dxdy\notag\\
& = O(1) +\epsilon^{-\frac{k}{k+1}}||U_1||^2_{S^2_{1, 0}(\mathbb R^2)}. \label{ueps_S}
\end{align}
From \eqref{S0} we deduce that $S_0\le ||u_\epsilon||^2_{S^2_{1, 0}(\Omega)}/||u_\epsilon||^2_{L^{p+1}_k(\Omega)}.$ On the other hand $U_1$ is an extremal function of Sobolev inequality \eqref{Sob} with the best constant $S$ we have
$$\lim_{\epsilon\to 0_+}\frac{O(1) +\epsilon^{-\frac{k}{k+1}}||U_1||^2_{S^2_{1, 0}(\mathbb R^2)}}{O(1)+\epsilon^{-\frac{k}{k+1}}||U_1||_{L^{p+1}_k(\mathbb R^2)}^2}=S.$$
Hence, from \eqref{ueps_p}-\eqref{ueps_S}, $S\le S_0.$ Using \eqref{S_0} we establish $S=S_0.$\\
Next we prove that $S_\lambda<S.$ Since $B_R\subset \Omega$ we have
\begin{align*}
||u_\epsilon||^2_{L^2_\beta(\Omega)} & \ge \int_{B_R} |x|^{2\beta}|\Phi_\epsilon(r)|^2dxdy\\
& \ge C_{k, \beta}\int_0^R r^{\frac{2\beta+1}{k+1}}(\epsilon^2+r^2)^{-\frac{k}{k+1}}dr. 
\end{align*}
Therefore we obtain the following statements:
\begin{itemize}
\item For $k=2(\beta+1)$ we have $||u_\epsilon||^2_{L^2_\beta(\Omega)}\ge C|\ln\epsilon|$.
\item For $k<2(\beta+1)$ we have  $||u_\epsilon||^2_{L^2_\beta(\Omega)}=O(1).$
\item For $k>2(\beta+1)$ we have $||u_\epsilon||^2_{L^2_\beta(\Omega)}\ge C\epsilon^{\frac{2(\beta+1)-k}{k+1}}.$
\end{itemize}
Besides, form \eqref{ueps_p}-\eqref{ueps_S} and \eqref{Slam} we get
\begin{align*}
S_\lambda & \le \frac{||u_\epsilon||^2_{S^2_{1, 0}(\Omega)} - \lambda||u_\epsilon||^2_{L^2_\beta(\Omega)}}{||u_\epsilon||^2_{L^{p+1}_k(\Omega)}}\\
& \le \frac{O(1)+\epsilon^{-\frac{k}{k+1}}||U_1||^2_{S^2_{1, 0}(\mathbb R^2)}-\lambda||u_\epsilon||^2_{L^2_\beta(\Omega)}}{O(1)+\epsilon^{-\frac{k}{k+1}}||U_1||^2_{L^{p+1}_k(\mathbb R^2)}}.
\end{align*}
Since $\lambda>0,$ $U_1$ is an extremal function of the Sobolev inequality \eqref{Sob} with the best constant $S=S_0$, for $k\ge 2(\beta+1)$ and let $\epsilon$ be small enough we are done.
\end{proof}
Now we can give the proof of part (c) of Theorem \ref{Case1}.
\begin{proof}[Proof of Theorem \ref{Case1} (c)]
From Remark \ref{Rem3}, we only need prove \eqref{Slam} has a nontrivial minimizer. There is a sequence $u_j\in S^2_{1, 0}(\Omega), j\in\mathbb N$ such that
\begin{equation}\label{mini}
||u_j||_{L^{p+1}_k(\Omega)}=1, ||u_j||^2_{S^2_{1, 0}(\Omega)}-\lambda||u_j||^2_{L^2_\beta(\Omega)}=S_\lambda+o(1).
\end{equation}
Case 1: For $(p+1)\beta\ge 2k$, using Proposition \ref{imb} we have the bounded embedding $L^2_\beta(\Omega)\hookrightarrow L^{p+1}_k(\Omega).$ So from \eqref{mini}  the sequence $\{u_j\}_{j\in\mathbb N}$ is bounded in $S^2_{1, 0}(\Omega).$\\
Case 2: For $(p+1)\beta<2k$  using the Holder inequality we deduce that
\begin{equation}\label{Hol}
\int_{|x|>\epsilon\atop (x, y)\in\Omega}|x|^{2\beta}|u_j|^2dxdy \le ||u_j||^2_{L^{p+1}_k(\Omega)}\epsilon^{2\beta-\frac{4k}{p+1}}|\Omega|^{\frac{p-1}{p+1}}.
\end{equation}
Since $\beta>-1/2$ there is a $\gamma\in (0, 1)$ such that $\gamma+2\beta>0.$ Let $\epsilon>0$. Using the Hardy inequality we have
\begin{equation}\label{Har1}
\int_{|x|<\epsilon\atop (x, y)\in\Omega}|x|^{2\beta}|u_j|^2dxdy \le \epsilon^{\gamma+2\beta}\int_{|x|<\epsilon\atop (x, y)\in\Omega} \frac{|u_j|^2}{|x|^\gamma}dxdy \le \epsilon^{\gamma+2\beta}||u_j||^2_{S^2_{1, 0}(\Omega)}.
\end{equation}
From \eqref{mini}-\eqref{Hol}-\eqref{Har1} we get
\begin{equation}
\frac{1}{\lambda}\left(||u_j||^2_{S^2_{1, 0}(\Omega)} -S_\lambda - o(1)\right)\le \epsilon^{\gamma+2\beta}||u_j||^2_{S^2_{1, 0}(\Omega)} + C_\epsilon|\Omega|^{\frac{p-1}{p+1}}.
\end{equation}
Since $\gamma+2\beta>0$, for $\epsilon>0$ small enough we obtain that the sequence $\{u_j\}_{j\in\mathbb N}$ is bounded in $S^2_{1, 0}(\Omega).$ 

From the boundedness of $\{u_j\}_{j\in\mathbb N}$ in $S^2_{1, 0}(\Omega)$,  using Proposition \ref{imb}, there is a subsequence of $\{u_j\}_{j\in\mathbb N}$ which is still denoted by $\{u_j\}_{j\in\mathbb N}$ and $u\in S^2_{1, 0}(\Omega)$ such that
\begin{align*}
u_j \rightharpoonup u & \text{ weakly in } S^2_{1, 0}(\Omega),\\
u_j \rightharpoonup u & \text{ weakly in } L^{p+1}_{k}(\Omega),\\
u_j \rightarrow u & \text{ strongly in } L^2_\beta(\Omega),\\
u_j \rightarrow u & \text{ a.e. in } \Omega.
\end{align*}
So $||u||_{L^{p+1}_k(\Omega)} \le \liminf_{j\to\infty}||u_j||_{L^{p+1}_k(\Omega)}=1.$ Using Brezis-Lieb lemma (\cite{Brezis:1983}) we have
$$||u_j||^{p+1}_{L^{p+1}_k(\Omega)}=||u||^{p+1}_{L^{p+1}_k(\Omega)}+||u-u_j||^{p+1}_{L^{p+1}_k(\Omega)}+o(1).$$
So from \eqref{Sob} and $p>1$ we deduce that
\begin{equation}\label{BL}
1 \le ||u||^2_{L^{p+1}_k(\Omega)} + \dfrac{1}{S}||u-u_j||^2_{S^2_{1, 0}(\Omega)}+o(1).
\end{equation}
Since $\lim_{j\to\infty}\langle\nabla_G(u_j-u), \nabla_Gu\rangle=0, \lim_{j\to\infty}||u_j||_{L^2_\beta(\Omega)}=||u||_{L^2_\beta(\Omega)}$, and using \eqref{mini} we have
$$||u_j-u||^2_{S^2_{1, 0}(\Omega)}+||u||^2_{S^2_{1, 0}(\Omega)} - \lambda||u||^2_{L^2_\beta(\Omega)}=S_\lambda+o(1).$$
Hence from \eqref{BL} we get
$$||u||^2_{S^2_{1, 0}(\Omega)} - \lambda ||u||^2_{L^2_\beta(\Omega)}\le S_\lambda ||u||_{L^{p+1}_k(\Omega)}^2+\left(\frac{S_\lambda}{S}-1\right)||u-u_j||^2_{S^2_{1, 0}(\Omega)}+o(1).$$
From \eqref{import}, let $j$ be large enough, we obtain
\begin{equation}\label{minSlam}
||u||^2_{S^2_{1, 0}(\Omega)}-\lambda||u||^2_{L^2_\beta(\Omega)}\le S_\lambda||u||^2_{L^{p+1}_k(\Omega)}.
\end{equation}
From \eqref{S0}, \eqref{import}, using \eqref{mini} again and let $j$ go to infinity we have
$$\lambda||u||^2_{L^2_\beta(\Omega)}\ge S_0-S_\lambda>0.$$ 
So $u\not\equiv 0.$ Hence, from \eqref{minSlam}, $U=u/||u||_{L^{p+1}_k(\Omega)}\not\equiv 0$ is a minimizer of \eqref{Slam}. The proof is done.
\end{proof}
\section{Proof of Theorem \ref{Case2}} 
From (h1)-(h2) we have some following properties of $h$:
\begin{itemize}
\item[(h3)] There is a $\mu_1>0$ such that
$$|h(x, y, \xi)|\le \xi^p+\mu_1\xi, \forall \xi>0 \text{ almost everywhere in } \Omega.$$
\item[(h4)] For each $\epsilon>0$ there is a $C=C(\epsilon)$ such that
$$|h(x, y, \xi)|\le C\xi^p+\epsilon \xi, \forall \xi>0 \text{ almost everywhere in } \Omega.$$
\item[(h5)] For each $\epsilon>0$ there is a $C=C(\epsilon)$ such that
$$|h(x, y, \xi)|\le \epsilon \xi^p+C, \forall \xi>0 \text{ almost everywhere in } \Omega.$$
\end{itemize}
So $H$ has following properties:
\begin{itemize}
\item[(H1)] For each $\epsilon>0$ there is a $C=C(\epsilon)$ such that
$$|H(x, y, \xi)|\le \frac{C}{p+1}\xi^{p+1}+\frac{\epsilon}{2} \xi^2, \forall \xi>0 \text{ almost everywhere in } \Omega.$$
\item[(H2)] For each $\epsilon>0$ there is a $C=C(\epsilon)$ such that
$$|H(x, y, \xi)|\le \frac{\epsilon}{p+1}\xi^{p+1}+C \xi, \forall \xi>0 \text{ almost everywhere in } \Omega.$$
\end{itemize}
For proving Theorem \ref{Case2} we need the following version of Mountain Pass Lemma.
\begin{Lem}(\cite{Brezis:1984})\label{MPL}
Let $\Phi$ be a $C^1$ functional on a Banach space $E$. Suppose there exists a neighborhood $\mathcal N$ of $0$ in $E$ and a constant $\rho$ such that 
\begin{itemize}
\item[(i)] $\Phi(u)\ge \rho$ for every $u$ in the boundary of $\mathcal N$;
\item[(ii)] $\Phi(0)=0$ and $\Phi(v)<\rho$ for some $v\not\in \mathcal N.$
\end{itemize}
Set
$$c=\inf_{P\in\mathcal P}\max_{w\in P}\Phi(w)\ge \rho$$
where 
$$\mathcal P=\{g\in C([0, 1]; E):\; g(0)=0, g(1)=v\}.$$
Then there is a sequence $\{u_j\}_{j\in\mathbb N}$ in $E$ such that
$$\lim_{j\to\infty}\Phi(u_j)=c, \lim_{j\to\infty}||\Phi'(u_j)||_{E^*}=0. $$
\end{Lem}
Next we use Lemma \ref{MPL} for $E=S^2_{1, 0}(\Omega)$ and the following functional
\begin{align}
\Phi(u)&=\frac{||u||^2_{S^2_{1, 0}(\Omega)}}{2}+\frac{\mu_1||u||^2_{L^2_k(\Omega)}}{2}-\frac{\mu_1||u^+||^2_{L^2_k(\Omega)}}{2}- \notag\\
& - \frac{||u^+||^{p+1}_{L^{p+1}_k(\Omega)}}{p+1}-\frac{\mu||u^+||_{L^{q+1}_\beta(\Omega)}^{q+1}}{q+1}-\int_\Omega|x|^{2k}H(x, y, u^+)dxdy \label{dePhi}
\end{align}
where $u^+=\max\{u, 0\},$ $\mu_1$ is the constant in (h3). It is easy to see that
$$\Phi(0)=0 \text{ and } \Phi(u)=\Psi(u), \forall u\in S^2_{1, 0}(\Omega), u\ge 0.$$
Firstly we will prove the above functional $\Phi$ satisfies the assumption of  \linebreak Lemma \ref{MPL}.
\begin{Lem}\label{PhiMPL}
Let $\Phi$ be defined as \eqref{dePhi}. Then $\Phi$ belongs to $\in C^1(S^2_{1, 0}(\Omega); \mathbb R)$ and satisfies (i)-(ii) in Lemma \ref{MPL}.
\end{Lem}
\begin{proof}
It is not difficult to see that $\Phi\in C^1(S^2_{1, 0}(\Omega); \mathbb R)$ and
\begin{align}
\langle \Phi'(u), v\rangle=\int_\Omega & \Big(\nabla_Gu\cdot\nabla_G v+\mu_1|x|^{2k}uv-\mu_1|x|^{2k}u^+v - |x|^{2k}(u^+)^pv - \notag\\
 - & \mu |x|^{2\beta}(u^+)^qv -|x|^{2k}h(x, y, u^+)v\Big)dxdy, \forall u, v\in S^2_{1, 0}(\Omega)\label{Phi'}.
\end{align}
From (H1) we obtain
\begin{equation}\label{Phig}
\Phi(u)\ge \frac{||u||^2_{S^2_{1, 0}(\Omega)}}{2}  -\frac{\epsilon||u^+||^2_{L^2_k(\Omega)}}{2}-\frac{(C+1)||u^+||^{p+1}_{L^{p+1}_k(\Omega)}}{p+1}-\frac{\mu||u^+||_{L^{q+1}_\beta(\Omega)}^{q+1}}{q+1}.
\end{equation}
Hence, using Sobolev embedding in Proposition \ref{imb}, the compact embedding $S^2_{1, 0}(\Omega)\hookrightarrow L^{q+1}_\beta(\Omega)$ and $u^+\le|u|$ , we get
\begin{equation}
\Phi(u)\ge \frac{(1-C_1\epsilon)||u||^2_{S^2_{1, 0}(\Omega)}}{2} -\frac{C_2(C+1)||u||^{p+1}_{S^2_{1, 0}(\Omega)}}{p+1}-\frac{C_3\mu||u||_{S^2_{1, 0}(\Omega)}^{q+1}}{q+1}.
\end{equation}
So let $\epsilon>0$ small enough, there exists $R>0, \rho>0$ such that
$$\Phi(u)>\rho, \forall ||u||_{S^2_{1, 0}(\Omega)}=R.$$
 This means that $\Phi$ satisfies (i) in Lemma \ref{MPL} for $\mathcal N=\{||u||_{S^2_{1, 0}(\Omega)}<R\}.$ Next we will prove that $\Phi$  satisfies (ii). It is easy to see $\Phi(0)=0$. Let $v_1\in S^2_{1, 0}(\Omega)$ such that $v_1\ge 0,$ $v_1\not\equiv 0$. We have
\begin{align}
\Phi(tv_1)=\Psi(tv_1) & = \frac{t^2}{2}||v_1||^2_{S^2_{1, 0}(\Omega)} - \frac{t^{p+1}}{p+1}||v_1||^{p+1}_{L^{p+1}_k(\Omega)}-\notag\\
& - \mu\frac{t^{q+1}}{q+1}||v_1||^{q+1}_{L^{q+1}_\beta(\Omega)}-\int_\Omega |x|^{2k}H(x, y, tv_1 )dxdy. 
\end{align}
So using (H2) it is not difficult to see that there $t_1>0$ large enough such that 
$$\Phi(v_0)<0, v_0=t_1v_1.$$
\end{proof}
We now prove Theorem \ref{Case2}.
\begin{proof}[Proof of Theorem \ref{Case2}]
From Lemma \ref{PhiMPL} and Lemma \ref{MPL} there is a sequence $\{u_j\}_{j\in\mathbb N}$ in $S^2_{1, 0}(\Omega)$ such that
$$\lim_{j\to\infty}\Phi(u_j)=c, \lim_{j\to\infty}||\Phi'(u_j)||_{(S^2_{1, 0}(\Omega))^*}=0$$
where
$$c=\inf_{P\in\mathcal P}\max_{w\in P}\Phi(w)\ge \rho$$
and
$$\mathcal P=\{g\in C([0, 1]; S^2_{1, 0}(\Omega)):\; g(0)=0, g(1)=v_0\}.$$
Note that $P_0=\{g:[0, 1]\to S^2_{1, 0}(\Omega):\; g(t)=tv_0\}\in \mathcal P$  and $\Phi(tv_0)=\Psi(tv_0),$ $\forall t\ge 0,$ from \eqref{asmpC2} we have
\begin{equation}\label{cS}
c < \frac{1+k}{2+3k}S^{\frac{2+3k}{2(1+k)}}.
\end{equation}
From $\lim_{j\to\infty} \Phi(u_j)=c$ we have
\begin{align}
\dfrac{||u_j||^2_{S^2_{1, 0}(\Omega)}}{2} & +\dfrac{\mu_1(||u_j||^2_{L^2_k(\Omega)}-||u^+_j||^2_{L^2_k(\Omega)})}{2} - \dfrac{||u^+||^{p+1}_{L^{p+1}_k(\Omega)}}{p+1}-\notag\\
& - \dfrac{\mu||u^+_j||^{q+1}_{L^{q+1}_\beta(\Omega)}}{q+1} - \int_\Omega|x|^{2k}H(x, y, u_j^+)dxdy=c+o(1). \label{Phi0}
\end{align}
From \eqref{Phi'} we have
\begin{align}
\langle \Phi'(u_j), u_j\rangle= & ||u_j||^2_{S^2_{1, 0}(\Omega)} +\mu_1(||u_j||^2_{L^2_k(\Omega)}-||u_j^+||^2_{L^2_k(\Omega)}) -||u_j^+||^{p+1}_{L^{p+1}_k(\Omega)}-\notag\\
- & \mu||u_j^+||^{q+1}_{L^{q+1}_\beta(\Omega)} -\int_\Omega |x|^{2k}h(x, y, u_j^+)dxdy. \label{Phi's}
\end{align}
From \eqref{Phi0}-\eqref{Phi's} we get
\begin{align}
\left(\frac{1}{2}-\frac{1}{p+1}\right)||u_j^+||^{p+1}_{L^{p+1}_k(\Omega)}  +\mu\left(\frac{1}{2}-\frac{1}{q+1}\right)||u_j^+||^{q+1}_{L^{q+1}_\beta(\Omega)}& - \notag\\
-\int_\Omega|x|^{2k}\left(H(x, y, u_j^+)-\frac{1}{2}h(x, y, u_j^+)\right)dxdy -  \frac{1}{2}\langle \Phi'(u_j), u_j\rangle &= c+o(1).\label{PhiS1}
\end{align}
Using (h3) and (H2) we have
\begin{align}
\int_\Omega|x|^{2k}|h(x, y, u_j^+)u_j^+|dxdy & \le \epsilon||u^+_j||^{p+1}_{L^{p+1}_k(\Omega)}+C||u_j^+||_{L^1_k(\Omega)},\label{guj}\\
\int_\Omega|x|^{2k}|H(x, y, u_j^+)|dxdy & \le \frac{\epsilon}{p+1}||u^+_j||^{p+1}_{L^{p+1}_k(\Omega)}+C||u_j^+||_{L^1_k(\Omega)}.\label{Guj}
\end{align}
So using \eqref{PhiS1} and Sobolev embedding in Proposition \ref{imb}, note that 
\begin{equation}\label{Phi'0}
|\langle \Phi'(u_j), u_j\rangle|\le ||\Phi'(u_j)||_{(S^2_{1, 0}(\Omega))^*}||u_j||_{S^2_{1, 0}(\Omega)}
\end{equation}
and $\lim_{j\to\infty}||\Phi'(u_j)||_{(S^2_{1, 0}(\Omega))^*}=0$, let $\epsilon$ be small enough we have
\begin{align}
||u^+_j||^{p+1}_{L^{p+1}_k(\Omega)}&\le C+C||u_j||_{S^2_{1, 0}(\Omega)}, \label{ujLp}\\
||u^+_j||^{q+1}_{L^{q+1}_\beta(\Omega)}&\le C+C||u_j||_{S^2_{1, 0}(\Omega)}.\label{ujLq}
\end{align}
Hence from \eqref{Phi0} we obtain
$$||u_j||^2_{S^2_{1, 0}(\Omega)}\le C+C||u_j||_{S^2_{1, 0}(\Omega)}.$$
Therefore $\{u_j\}_{j\in\mathbb N}$ is bounded in $S^2_{1, 0}(\Omega).$ From continuous embeddings 
$$S^2_{1, 0}(\Omega)\hookrightarrow L^{p+1}_k(\Omega), S^2_{1, 0}(\Omega)\hookrightarrow L^{q+1}_\beta(\Omega)$$
$\{u_j\}_{j\in\mathbb N}$ is also bounded in $L^{p+1}_k(\Omega)$ and $L^{q+1}_\beta(\Omega)$. Moreover, using compact embedding $S^2_{1, 0}(\Omega)\hookrightarrow L^r_k(\Omega), 1\le r<p+1,$ we can extract a subsequence, still denoted by $\{u_j\}_{j\in\mathbb N}$, and there is $u\in S^2_{1, 0}$ so that
\begin{align*}
u_j \rightharpoonup u & \text{ weakly in } S^2_{1, 0}(\Omega),\\
u_j \rightarrow u & \text{ strongly in } L^{r}_{k}(\Omega),\\
u_j \rightarrow u & \text{ a.e. in } \Omega,\\
h(\cdot, \cdot, u_j^+) \rightarrow h(\cdot, \cdot, u^+) & \text{ a.e. in } \Omega.
\end{align*} 
So from (h5), \eqref{ujLp}-\eqref{ujLq} and using Lebesgue dominated convergence theorem we have
\begin{align*}
\lim_{j\to\infty} \int_\Omega|x|^{2k}(u_j^+)^p\varphi dxdy & = \int_\Omega|x|^{2k}(u^+)^p\varphi dxdy, \forall \varphi \in S^2_{1, 0}(\Omega)\hookrightarrow L^{p+1}_k(\Omega),\\
\lim_{j\to\infty}\int_\Omega|x|^{2\beta}(u_j^+)^q\varphi dxdy & =\int_\Omega|x|^{2\beta}(u_j^+)^q\varphi dxdy, \forall \varphi\in S^2_{1, 0}(\Omega)\hookrightarrow L^{q+1}_\beta(\Omega,\\
\lim_{j\to\infty}\int_\Omega|x|^{2k}h(x, y, u_j^+)\varphi dxdy & = \int_\Omega|x|^{2k}h(x, y, u_j^+)\varphi dxdy, \forall \varphi \in S^2_{1, 0}(\Omega)\hookrightarrow L^{p+1}_k(\Omega).
\end{align*}
Hence let $j\to\infty$ in \eqref{Phi'} with $u=u_j, v=\varphi$ we obtain
\begin{equation}\label{PTdc}
-\Delta_G u + \mu_1|x|^{2k} u= |x|^{2k}(\mu_1u^++(u^+)^p+h(x, y, u^+)) +\mu|x|^{2\beta}(u^+)^q \text{ in } (S^2_{1, 0}(\Omega))^*. 
\end{equation} 
From (h3) the right-hand side of \eqref{PTdc} is greater or equal to $0$. So using Maximum Principle (\cite{Monti:2009}) we have $u\ge 0$ a.e. in $\Omega.$ Hence $u$ satisfies \eqref{equ1} and \eqref{equ3}.

In order to prove that $u$ satisfies \eqref{equ2}, using Strong Maximum Principle we will only prove that $u$ is nontrivial. Indeed, suppose that $u\equiv 0$. We have 
\begin{align*}
u_j \rightharpoonup 0 & \text{ weakly in } S^2_{1, 0}(\Omega),\\
u_j \rightarrow 0 & \text{ strongly in } L^{r}_{k}(\Omega),1 \le r<p+1.
\end{align*}
Note that $\lim_{j\to\infty}||\Phi'(u_j)||_{(S^2_{1, 0}(\Omega))^*}=0$, from \eqref{guj}-\eqref{Guj}-\eqref{Phi'0} and the compact embedding $S^2_{1, 0}(\Omega)\hookrightarrow L^{q+1}_\beta(\Omega)$ we have
\begin{align*}
\lim_{j\to\infty} \langle \Phi'(u_j), u_j\rangle & =0,\\
\lim_{j\to\infty} \int_\Omega |x|^{2k}|h(x, y, u_j^+)u_j^+|dxdy & =0,\\
\lim_{j\to\infty} \int_\Omega |x|^{2k}|H(x, y, u_j^+)|dxdy & =0,\\
\lim_{j\to\infty} \int_\Omega |x|^{2\beta}(u_j^+)^{q+1}dxdy & =0.
\end{align*}
So let $j\to\infty$ in \eqref{Phi's} we have
$$\lim_{j\to\infty} \left(||u_j||^2_{S^2_{1, 0}(\Omega)} -||u_j^+||^{p+1}_{L^{p+1}_k(\Omega)}\right)=0.$$
Hence there is an $\ell>0$ such that
\begin{equation}\label{L1}
\lim_{j\to\infty} ||u_j||^2_{S^2_{1, 0}(\Omega)} =\lim_{j\to\infty}||u_j^+||^{p+1}_{L^{p+1}_k(\Omega)}=\ell.
\end{equation}
Then let $j\to\infty$ in \eqref{Phi0} we get
\begin{equation}
\left(\frac{1}{2}-\frac{1}{p+1}\right)\ell =c.
\end{equation}
From Sobolev inequality \eqref{Sob} we have $S||u_j||^2_{L^{p+1}_k(\Omega)}\le ||u_j||^2_{S^2_{1, 0}(\Omega)}.$ Let $j\to \infty$ we deduce that
\begin{equation}\label{L2}
\ell^{2(\frac{1}{2}-\frac{1}{p+1})}\ge S.
\end{equation}
Since $p=(4+5k)/k$, from \eqref{L1}-\eqref{L2} we obtain
$$c \ge \frac{1+k}{2+3k}S^{\frac{2+3k}{2(1+k)}}$$
which contradicts \eqref{cS}. Therefore the proof is done.
\end{proof}
For $k>0, \beta>-1/2, 1<q<p,$ and
$$2k(q-1)+q<2\beta(p-1)+p,$$
from Proposition \ref{imb} we have the compact embedding $S^2_{1, 0}(\Omega)\hookrightarrow L^{q+1}_\beta(\Omega).$ So for proving Corollary \ref{Cor} we will show that \eqref{asmpC2} holds as follows. Firstly we prove Corollary \ref{Cor} (a).
\begin{proof}[Proof of Corollary \ref{Cor} (a)]
Let $v_\epsilon=u_\epsilon/||u_\epsilon||_{L^{p+1}_k(\Omega)},$ where $u_\epsilon$ is defined in \linebreak Section 4 (Proof of Theorem \ref{Case1}). From \eqref{ueps_p}-\eqref{ueps_S} we have 
\begin{equation}\label{veps}
||v_\epsilon||_{L^{p+1}_k(\Omega)}=1 \text{  and } ||v_\epsilon||_{S^2_{1, 0}(\Omega)}^2=S+O(\epsilon^{\frac{k}{k+1}}).
\end{equation}
Since $\Omega$ is bounded, there is an $R_2>0$ such that $\Omega\subset B_{R_2}.$ Then
\begin{align*}
||u_\epsilon||^{q+1}_{L^{q+1}_\beta(\Omega)} &\ge \int_{B_{R_2}} |x|^{2\beta}|\Phi_\epsilon|^{q+1}dxdy\\
& \ge C_{k, \beta} \int_0^{R_2}r^{\frac{2\beta+1}{k+1}}(\epsilon^2+r^2)^{-\frac{(q+1)k}{2(k+1)}}dr.
\end{align*} 
So, note that $4=(p-1)k/(k+1),$ we have 
\begin{equation}\label{veps1}
||v_\epsilon||^{q+1}_{L^{q+1}_\beta(\Omega)} \le 
\begin{cases}C\epsilon^{\frac{(q+1)k}{2(k+1)}}|\ln\epsilon| & \text{ for } 2\beta+2=qk,\\ 
C\left(\epsilon^{\frac{(q+1)k}{2(k+1)}}+\epsilon^{\frac{k[(p-1)(\beta+1)-(q-1)(k+1)]}{2(k+1)^2}}\right) & \text{ for } 2\beta+2\not=qk.\end{cases}
\end{equation}
Since $2k(q-1)+q<2\beta(p-1)+p, 1<q<p$ and $k>0$ we have
\begin{equation}\label{veps2}
\lim_{\epsilon\to 0_+} ||v_\epsilon||_{L^{q+1}_\beta(\Omega)}=0.
\end{equation}
From \eqref{veps} and $h=0$ we get
$$\Psi(tv_\epsilon)=\frac{t^2}{2}X_\epsilon - \frac{t^{p+1}}{p+1}-\mu\frac{t^{q+1}}{q+1}||v_\epsilon||^{q+1}_{L^{q+1}_\beta(\Omega)}$$
where $X_\epsilon=||v_\epsilon||^2_{S^2_{1, 0}(\Omega)}=S+O(\epsilon^{\frac{k}{k+1}}).$ In view of
$$\lim_{t\to+\infty}\Psi(tv_\epsilon)=-\infty, \Psi(0)=0,$$
and $\Psi(tv_\epsilon)>0$ for some $t>0,$
there exists $t_\epsilon>0$ such that
$$\Psi(t_\epsilon v_\epsilon)=\sup_{t\ge 0}\Psi(tv_\epsilon)=:Y_\epsilon.$$
Then $t_\epsilon X_\epsilon-t_\epsilon^p=\mu t_\epsilon^q||v_\epsilon||^{q+1}_{L^{q+1}_\beta(\Omega)}$ and $0<t_\epsilon<X_\epsilon^{1/(p-1)}.$ From \eqref{veps}-\eqref{veps2} and $q>1$ we have $\lim_{\epsilon\to 0_+}t_\epsilon=S^{1/(p-1)}.$ Moreover, note that $p=(4+5k)/k$ we obtain
\begin{align}
Y_\epsilon & = \frac{t_\epsilon^2}{2}X_\epsilon -\frac{t_\epsilon^{p+1}}{p+1}-\mu\frac{t_\epsilon^{q+1}}{q+1}||v_\epsilon||^{q+1}_{L^{q+1}_\beta(\Omega)}\notag\\
& \le \left(\frac{1}{2}-\frac{1}{p+1}\right)X_\epsilon^{\frac{p+1}{p-1}}-\mu\frac{t_\epsilon^{q+1}}{q+1}||v_\epsilon||^{q+1}_{L^{q+1}_\beta(\Omega)}\notag\\
&\le \frac{1+k}{2+3k}S^{\frac{2+3k}{2(1+k)}}+O(\epsilon^{\frac{k}{k+1}})-\mu\frac{t_\epsilon^{q+1}}{q+1}||v_\epsilon||^{q+1}_{L^{q+1}_\beta(\Omega)}\label{Yeps}.
\end{align}
Note that $B_R\subset\Omega$, by the same way as \eqref{veps1} we have
\begin{equation}\label{veps11}
\epsilon^{-\frac{k}{k+1}}||v_\epsilon||^{q+1}_{L^{q+1}_\beta(\Omega)} \ge 
\begin{cases}\bar{C}\epsilon^{\frac{(q-1)k}{2(k+1)}}|\ln\epsilon| & \text{ for } 2\beta+2=qk,\\ 
\bar{C}\left(\epsilon^{\frac{(q-1)k}{2(k+1)}}+\epsilon^{\frac{4(\beta+1)-(q+1)k}{2(k+1)}}\right) & \text{ for } 2\beta+2\not=qk.\end{cases}
\end{equation}
Since $q>1, k>0$ and $(q+1)k>4(\beta+1)$ then $qk\not=2\beta+2.$ Hence, note that $\mu>0,$ there is an $\epsilon>0$ such that
$$Y_\epsilon< \frac{1+k}{2+3k}S^{\frac{2+3k}{2(1+k)}} \text{ or } \sup_{t>0}\Psi(tv_\epsilon)<\frac{1+k}{2+3k}S^{\frac{2+3k}{2(1+k)}}.$$ 
Therefore \eqref{asmpC2} is proved.
\end{proof}
Next we prove Corollary \ref{Cor} (b).
\begin{proof}[Proof of Corollary \ref{Cor} (b)]
Let $v_1=u_1/||u_1||_{L^{p+1}_k(\Omega)},$ where $u_\epsilon, \epsilon=1,$ is defined in Section 4 (Proof of Theorem \ref{Case1}). We have
\begin{itemize}
\item $||v_1||_{L^{p+1}_k(\Omega)}=1, v_1\ge 0$.
\item $\Psi_\mu(tv_1)=\Psi(tv_1)=\frac{t^2}{2}||v_1||^2_{S^2_{1, 0}(\Omega)} -\frac{t^{p+1}}{p+1}-\mu\frac{t^{q+1}}{q+1}||v_1||^{q+1}_{L^{q+1}_\beta(\Omega)}.$
\end{itemize}
For each $\mu>0$ we have 
$$\Psi_\mu(0)=0, \lim_{t\to+\infty}\Psi_\mu(tv_1)=-\infty,$$
and $\Psi_\mu(tv_1)>0$ for some $t>0,$
so there is a $t_\mu>0$ such that $\Psi_\mu(t_\mu v_1)=\sup_{t>0}\Psi_\mu(tv_1).$ Then
$$t_\mu ||v_1||^2_{S^2_{1, 0}(\Omega)} - t_\mu^p = \mu t_\mu^q||v_1||^{q+1}_{L^{q+1}_\beta(\Omega)}, \text{ and } 0<t_\mu \le ||v_1||^{\frac{2}{p-1}}_{S^2_{1, 0}(\Omega)}.$$
Hence $\lim_{\mu\to+\infty}t_\mu=0$ and
$$\limsup_{\mu\to+\infty}\sup_{t>0} \Psi_\mu(tv_1)\le \lim_{\mu\to+\infty}\left(\frac{t_\mu^2}{2}||v_1||^2_{S^2_{1, 0}(\Omega)} -\frac{t_\mu^{p+1}}{p+1}\right)=0.$$
Therefore \eqref{asmpC2} holds for $\mu>\mu_0,$ where $\mu_0$ is enough large. The proof is done.
\end{proof}
\section{Remarks on nonexistence solution for Case 2}
Let us now consider the problem \eqref{equ1}-\eqref{equ3} for Case 2 and $h=0$, $\Omega$ is G-starshaped. By the same way to have \eqref{Pohoz} we have another Pohozaev-type identity as follows.
\begin{Md}
Assume that $u\in S_2(\overline{\Omega})$ is a solution of the problem \eqref{equ1}-\eqref{equ3}. Then we have
\begin{equation}\label{Pohoz1}
\frac{1}{2}\int_{\partial\Omega}(\mathcal T\cdot\nu)(\nu_1^2+|x|^2\nu_2^2)|\partial_\nu u|^2dS=\frac{\mu[4(\beta+1)-(q-1)k]}{2(q+1)}\int_\Omega|x|^{2\beta}u^{q+1}dxdy.
\end{equation}
\end{Md}
\noindent As in \cite{Lanconelli:2012}, using Pohozaev-type identity \eqref{Pohoz1} we get the following nonexistence solution result for Case 2.
\begin{Th}\label{DL_non}
Assume that $\Omega$ is G-starshaped, $h=0$, $\beta>-1/2,$ $1< q<p$ and
\begin{equation}\label{non_e}
\mu[4(\beta+1)-(q-1)k]\le 0.
\end{equation}
Then the problem \eqref{equ1}-\eqref{equ3} does not have a solution in $S_2(\overline{\Omega}).$
\end{Th}
\begin{Cor}\label{hq22}
Let $\Omega$ be G-starshaped, $h=0,$ $\beta>-1/2$ and $1< q<p$ such that 
$$2k(q-1)+q<2\beta(p-1)+p.$$
\begin{itemize}
\item For $\mu>0,$ $4(\beta+1)<(q+1)k$, or $\mu>\mu_0,$ with $\mu_0>0$ large enough, the problem \eqref{equ1}-\eqref{equ3} has a solution in $S^2_{1, 0}(\Omega).$ \item For $\mu\le 0$ the  problem \eqref{equ1}-\eqref{equ3} has no solution in $S_2(\overline{\Omega}).$
\end{itemize}
\end{Cor}
\begin{proof}
 Note that $4=(p-1)k/(k+1),$ for $2k(q-1)+q<2\beta(p-1)+p$ we have
$$4(\beta+1)-(q-1)k>0.$$
Hence from Corollary \ref{Cor} and Theorem \ref{DL_non} the conclusion follows.
\end{proof}
Next we consider the case $\mu>0, 4(\beta+1)>(q+1)k$ and $\Omega$ is strictly G-starshaped, i.e. there exists an $\epsilon_0>0$ such that
$$\mathcal T\cdot \nu \ge \epsilon_0 \text{ on } \partial\Omega.$$
Assume that $u\in C^2(\overline{\Omega})$ is a solution of the problem \eqref{equ1}-\eqref{equ3}. Using Divergence Theorem we have
$$\int_\Omega \Delta_Gudxdy =\int_{\partial\Omega}(\nu_1^2+|x|^{2k}\nu_2^2)\partial_\nu udS.$$
Note that $\mathcal T\cdot\nu\ge\epsilon_0$ on $\partial\Omega,$ $(-\Delta_Gu)>0$ in $\Omega$, and using Pohozaev-type identity \eqref{Pohoz1} and the Holder inequality we have
\begin{equation}\label{nonE1}
\left(\int_\Omega |\Delta_Gu|dxdy\right)^2\le \dfrac{C\mu[4(\beta+1)-(q-1)k]}{q+1}\int_\Omega|x|^{2\beta}u^{q+1}dxdy
\end{equation}
where $C=\int_{\partial\Omega}(\nu_1^2+|x|^{2k}\nu_2^2)/(\mathcal T\cdot\nu)dS.$ Here, for $k\in\mathbb N,$ we have the fundamental solution $F_{k, 0}(x, y, z, t)$ of $\Delta_G$ (\cite{N. M. Tri:1998a}, \cite{N. M. Tri:2002}) and some properties of it as follows.
\begin{itemize}
\item[(P1)] $\Delta_G F_{k, 0}=\delta(x-z, y-t).$
\item[(P2)] For each compact set $K\subset \mathbb R^2$, there is an $r_0>0$ such that
$$F_{k, 0}(x, y, z, t)\le 0, \forall (z, t)\in K, \forall (x, y)\in B_{r_0}(z, t),$$
where $B_{r_0}(z, t)=\{(x, y)\in\mathbb R^2:\; (x-z)^{2k+2}+(k+1)^2(y-t)^2<r_0^2\}.$
\item[(P3)] For every $(x, y, z, t)\in\mathbb R^4\setminus\{(x, y)=(z, t)\}$ we have
$$|F_{k, 0}(x, y, z, t)|\le C((x-z)^{2k+2}+(k+1)^2(y-t)^2)^{-\frac{k}{2k+2}}.$$
\end{itemize}
Let $K=\overline{\Omega},$ and $r_0>0$ as in (P2). It is obvious to see that $g: \overline{\Omega}\to \mathbb R$ defined by
$$g(z, t)=\max_{(x, y)\in\overline{\Omega}\setminus B_{r_0}(z, t)}F_{k, 0}(x, y, z, t)$$
is continuous. So there is a $C_0>0$ such that
\begin{equation}\label{fund_Tri}
F_{k, 0}(x, y, z, t)\le C_0, \forall (x, y), (z, t)\in\overline{\Omega}, (x, y)\not=(z, t).
\end{equation}
Put $\mathcal K(x, y, z, t)=C_0-F_{k, 0}(x, y, z, t)$. From (P1), \eqref{fund_Tri} and Maximum Principle we get
\begin{equation}\label{Fund}
u(x, y) \le \int_\Omega \mathcal K(x, y, z, t)|\Delta_Gu(z, t)|dzdt, \forall (x, y)\in\Omega.
\end{equation}
Since $\Omega$ is bounded, using (P3) we have
$$\int_{(x, y)\in\Omega\atop \mathcal K(x, y, z, t)>\lambda}|x|^{2k}dxdy \le C\int_0^{\lambda^{-\frac{k+1}{k}}}r^{\frac{1}{k+1}}dr=C\lambda^{-\frac{p-3}{2}}.$$
Hence 
$$\sup_{(z, t)\in\Omega} ||\mathcal K(x, y, z, t)||_{L^{\frac{p-3}{2}, \infty}_k(\Omega)}<\infty,$$
where $||w||_{L^{\frac{p-3}{2}, \infty}_k(\Omega)}=\sup_{\lambda>0}\lambda\left(\int_{(x, y)\in\Omega\atop|w|>\lambda}|x|^{2k}dxdy\right)^{2/(p-3)}.$ So from \eqref{Fund} we have
\begin{equation}\label{weak1}
||u||_{L^{\frac{p-3}{2}, \infty}_k(\Omega)}\le C\int_\Omega|\Delta_Gu(z, t)|dzdt.
\end{equation}
We need the following technique lemma.
\begin{Lem}\label{teq}
Assume that $\beta>-1/2,$ $1\le q<p-2,$ and 
\begin{equation}\label{nasi2}
(p-1)(\beta+1/2)>(q+1)(k+1/2).
\end{equation}
Then there is a $q_0\in [q, p-2)$  such that
\begin{equation}\label{ns0}
\int_\Omega |x|^{2\beta}|u|^{q+1}dxdy \le C_s\left(\int_\Omega|x|^{2k}|u|^{p}dxdy\right)^{\frac{\theta(q+1)}{p}}||u||_{L^{s, \infty}_k(\Omega)}^{(1-\theta)(q+1)}, s\in[1, q_0+1),
\end{equation}
where $\theta=\frac{(q_0+1-s)p}{(p-s)(q_0+1)}.$
\end{Lem}
\begin{proof}
For $k\le \beta$, since $\Omega$ is bounded, let $q_0=q$ we have
\begin{equation}\label{ns1}
\int_\Omega |x|^{2\beta}|u|^{q+1}dxdy\le C\int_\Omega |x|^{2k}|u|^{q_0+1}dxdy.
\end{equation}
For $k>\beta$, since $\beta>-1/2, 1\le q<p-2$ and $(p-1)(\beta+1/2)>(q+1)(k+1/2)$, there is a $q_0\in(q, p-2)$ such that
\begin{equation}\label{nasi1}
(q_0+1)(\beta+1/2)>(q+1)(k+1/2).
\end{equation}
So from Proposition \ref{imb} we have
\begin{equation}\label{ns2}
\int_\Omega |x|^{2\beta}|u|^{q+1}dxdy\le C\left(\int_\Omega |x|^{2k}|u|^{q_0+1}dxdy\right)^{\frac{q+1}{q_0+1}}.
\end{equation}
Since $1\le s<q_0+1<p$ we have the interpolation inequality (see \cite{Graf})
$$\int_\Omega |x|^{2k}|u|^{q_0+1}dxdy \le \frac{(q_0+1)(p-s)}{(q_0+1-s)(p-q_0-1)}||u||_{L^{p, \infty}_k(\Omega)}^{\theta(q_0+1)}||u||_{L^{s, \infty}_k(\Omega)}^{(1-\theta)(q_0+1)}.$$
Hence from \eqref{ns1}-\eqref{ns2} and $||u||_{L^{p, \infty}_k(\Omega)}\le C||u||_{L^p_k(\Omega)}$ we get \eqref{ns0}.
\end{proof}
\begin{Md}\label{non_e2}
Assume that $\Omega$ is strictly G-starshaped, $\beta>-1/2,$ $1< q<p-2$ and \eqref{nasi2} holds.
If there exists a $q_0\in [q, p-2)$ and an $s\in [1, (p-3)/2)$ such that
\begin{equation}\label{non_e3}
\frac{(q+1)[q_0+1+(p-q_0-2)s]}{(q_0+1)(p-s)}=2
\end{equation}
and Lemma \ref{teq} holds
then there exists a $\mu_1>0$ such that the problem \eqref{equ1}-\eqref{equ3} does not have a solution in $C^2(\overline{\Omega})$ for $0<\mu<\mu_1.$
\end{Md}
 \begin{proof}
Since $4=(p-1)k/(k+1)$, if $(p-1)(\beta+1/2)>(q+1)(k+1/2)$  and $p-q>2$ then
$$4(\beta+1)>(q+1)k>(q-1)k.$$
Assume that for some $\mu>0$ the problem \eqref{equ1}-\eqref{equ3} has a solution $u\in C^2(\overline{\Omega}).$ Then we have \eqref{nonE1}-\eqref{weak1} and
$$|\Delta_Gu|\ge |x|^{2k}u^p \text{ in } \Omega.$$
From \eqref{non_e3} and $p-2>q_0\ge q> 1$ we have 
$$1\le s=\dfrac{\frac{2(q_0+1)p}{q+1}-(q_0+1)}{\frac{2(q_0+1)}{q+1}+p-(q_0+2)}<q_0+1.$$
Hence from Lemma \ref{teq} we have \eqref{ns0}. For $s\le (p-3)/2$ we have
$$||u||_{L^{s, \infty}_k(\Omega)}\le C||u||_{L^{(p-3)/2, \infty}_k(\Omega)}.$$
Therefore 
\begin{align*}
\int_\Omega |x|^{2\beta}u^{q+1}dxdy & \le C\left(\int_\Omega|x|^{2k}u^pdxdy\right)^{\frac{(q_0+1-s)(q+1)}{(p-s)(q_0+1)}}||u||_{L^{s, \infty}_k(\Omega)}^{\frac{[p-(q_0+1)](q+1)s}{(p-s)(q_0+1)}} \; (\text{from \eqref{ns0}})\\
& \le C\left(\int_\Omega|\Delta_Gu|dxdy\right)^{\frac{(q_0+1-s)(q+1)}{(p-s)(q_0+1)}}||u||_{L^{(p-3)/2, \infty}_k(\Omega)}^{\frac{[p-(q_0+1)](q+1)s}{(p-s)(q_0+1)}} \\
& \le C(\int_\Omega |\Delta_Gu|dxdy)^2 \; (\text{from \eqref{weak1}-\eqref{non_e3}})\\
& \le C\mu\int_\Omega |x|^{2\beta}u^{q+1}dxdy \; (\text{ from \eqref{nonE1}}).
\end{align*}
Hence  $C\mu \ge 1$, the proof is done.
\end{proof}
\begin{Cor}
 Let $\Omega$ be strictly G-starshaped, $\beta>-1/2$, $1< q<p-2$ and one of the following conditions holds.
\begin{itemize}
\item[(i)] $k\le \beta$ and $q\le p-2(p-1)/(p-5).$
\item[(ii)] $k>\beta$ and 
\begin{equation}\label{cond_n}
\frac{k+1/2}{\beta+1/2}<\frac{(p-1)(p-3)}{2(p+3)+(q+1)(p-5)}.
\end{equation}
\end{itemize}
Then there exists $\mu_0, \mu_1$ such that $0<\mu_1<\mu_0$ and
\begin{itemize}
\item for $\mu\le 0$  the problem \eqref{equ1}-\eqref{equ3} has no solution $u\in S_2(\overline{\Omega});$
\item for $0<\mu<\mu_1$ the problem \eqref{equ1}-\eqref{equ3} has no solution $u\in C^2(\overline{\Omega});$
\item for $\mu>\mu_0$ the problem \eqref{equ1}-\eqref{equ3} has a solution $u\in S^2_{1, 0}({\Omega}).$
\end{itemize}
\end{Cor}
\begin{proof}
For $k\le \beta$ we take $q_0=q$. Then 
$$s=\frac{2p-q-1}{p-q}\le \frac{p-3}{2} \Longleftrightarrow q\le p-\frac{2(p-1)}{p-5}.$$
Hence 
\begin{equation}\label{ptc}
2k(q-1)+q<2\beta(p-1)+p \text{ and } 4(\beta+1)>(q+1)k>(q-1)k.
\end{equation}
For $k>\beta$, in order to have a $q_0\in [q, p-2)$ and an $s\in [1, (p-3)/2)$ such that \eqref{nasi1} and \eqref{non_e3} hold we need \eqref{cond_n}.
Note that from \eqref{cond_n} we have \eqref{nasi2}. Then
\eqref{ptc} holds.

From \eqref{ptc}, Corollary \ref{hq22} and Proposition \ref{non_e2} we get the conclusion.
\end{proof}

\end{document}